\documentclass[12pt,leqno]{article}
\usepackage{amsmath,amsthm}

\def\init{\setcounter{equation}{0}}
\newtheorem{theorem}{Theorem}[section]

\newcommand{\R}{{\bf R}}

\newtheorem{pro}{Proposition}[section]

\newcommand{\e}{{\varepsilon}}

\title{Artificial black holes.
\author{G.Eskin, \ \ \  Department of Mathematics, UCLA,\\ Los Angeles,
CA 90095-1555, USA. \ E-mail: eskin@math.ucla.edu}
}

\begin{document}

\maketitle
\ \ \ \ \ \ \ \ \ \ \ \ \ \ \ \ \ \ \ \ \ \ \ \ \ \ 
To Misha Shubin on the occasion of his 65th birthday.
\\ 
\ 
\\
\begin{abstract}
We study black holes for the linear hyperbolic equations describing the wave propagation in
the moving medium.  Such black holes are called artificial since the Lorentz metric associated with 
the hyperbolic equation does not necessary satisfies the Einstein equations.  Artificial black holes 
also arise when we consider perturbations of the Einstein equations. In this paper we review results of 
[E2] and [E3] on the existence and the stability of black holes for the stationary wave equations in two
space dimensions, and in the axisymmetric case.  
\end{abstract}

\section{Introduction.}
\label{section 1}
\init

Consider the wave equation of the form
\begin{equation}                                               \label{eq:1.1}
\sum_{j,k=0}^n\frac{1}{\sqrt{(-1)^ng(x)}}\frac{\partial}{\partial x_j}\left(\sqrt{(-1)^ng(x)}g^{jk}(x)
\frac{\partial u(x_0,x)}{\partial x_k}\right)=0,
\end{equation}
where $x=(x_1,...,x_n)\in \R^n,\ x_0\in\R$  is the time variable,
the coefficients $g^{jk}(x)\in C^\infty$  and are independent of $x_0$,
 $g(x)=\det[g_{jk}(x)]_{j,k=0}^n,\ [g_{jk}(x)]_{j,k=0}^n$ 
is the inverse to the matrix $[g^{jk}(x)]_{j,k=0}^n$.
We assume  that  $[g_{jk}(x)]_{j,k=0}^n$ is a pseudo-Riemanian metric  with the Lorentz signature 
$(1,-1,...,-1)$.
We also assume that 
\begin{equation}                                               \label{eq:1.2}
g^{00}(x)>0,\ \ \forall x\in\R^n
\end{equation}
and
\begin{equation}                                               \label{eq:1.3}
g^{jk}(x)-\delta_{jk}=O\left(\frac{1}{|x|}\right)\ \ \mbox{when}\ \ \ |x|\rightarrow \infty.
\end{equation}
Equation (\ref{eq:1.1})  describes the wave propagation in a moving medium.
 As in [E1], [E2],  we consider two examples:
\\
\ 
\\
a) \underline{Propagation of light in a moving dielectric medium (cf. [G], [LP]).}
\\
\\
  In this case equation (\ref{eq:1.1}) is called the Gordon equation and it has the form:
\begin{equation}                                   \label{eq:1.4}
g^{jk}(x)=\eta^{jk}+(n^{2}(x)-1)u^ju^k,\ \ 0\leq j,k\leq n,\ \ n=3,
\end{equation}
 $\eta^{jk}=0$ when $j\neq k,\ \eta^{00}=1,\ \eta^{jj}=-1$ when $ 1\leq j\leq 3, 
\ n(x)=\sqrt{\e(x)\mu(x)}$  is  the refraction index,
$x_0=ct$  where $c$ is the speed of light in the vacuum, $(u_0,u_1,u_2,u_3)$  is the four-velocity:
$u_0=\left(1-\frac{|w|^2}{c^2}\right)^{-\frac{1}{2}},\  
u_j(x)=\frac{w_j(x)}{c}\left(1-\frac{|w_j|^2}{c^2}\right)^{-\frac{1}{2}},
\ 1\leq j\leq 3$.  Here   $w=(w_1(x),w_2(x),w_3(x))$  is the velocity of the flow, 
 $|w|=\sqrt{w_1^2+w_2^2+w_3^2}$.
\\
\
\\
b) \underline{Acoustic waves in a moving fluid  
(cf. [V1]).}
\\
\\
 In this case 

\begin{eqnarray}                                     \label{eq:1.5}
g^{00}=\frac{1}{\rho c},\ \ g^{0j}=g^{j0}=\frac{1}{\rho c}v^j,\ 1\leq j\leq 3,
\\
\nonumber
g^{jk}=-\frac{1}{\rho c}(-c^2\delta_{ij}+v^jv^k),\ \ 1\leq j,k\leq 3,
\end{eqnarray}
where $\rho$  is the density,  $c$  is the sound speed, $v=(v^1,v^2,v^3)$  is the velocity.  

Equations with metrics  (\ref{eq:1.4}), (\ref{eq:1.5})  may have black holes (see \S 2 below).  
These black holes
are called optical and acoustic black holes,  respectively (cf. [V1], [NVV], [U] 
and references there).  They are called often  artificial black holes since the metric in (\ref{eq:1.1})
not necessarily satisfies the Einstein equations.
Physicists hope to create artificial black holes in the laboratory. The artificial black holes play role 
when one consider the perturbations of black holes of the general relativity such as the 
Schwarzschield and the Kerr black holes.  We introduce the black holes in \S 2 and \S 3.  In \S 4
we shall study the existence and the stability of black holes in the case of two space dimensions.
In \S 5  we consider the axisymmetric case and in \S 6  we consider the inverse problems in
the presence of black holes.

\section{The black and  white holes.}
\label{section 2}
\init

Let $S_0(x)=0$  be a closed bounded smooth surface in $\R^n$.  Let $\Omega_{int}$  and $\Omega_{ext}$
be the interior and the exterior of $S_0(x)=0$,  respectively.
The domain $\Omega_{int}\times \R$ is called a black hole for  (\ref{eq:1.1})  if no signals (disturbances) 
from $\Omega_{int}\times \R$  can reach the exterior domain $\Omega_{ext}\times\R$.  Analogously, 
$\Omega_{int}\times \R$  is a white hole for (\ref{eq:1.1}) if no signals (disturbances)  from 
$\Omega_{ext}\times\R$
can reach $\Omega_{int}\times\R$.  The surface $\{S_0(x)=0\}\times \R$  is called the event horizon
if $\Omega_{int}\times\R$  is either black or white hole.
In order to find the conditions when $\{S_0(x)=0\}\times \R$  is an event horizon we need the notion of
the forward domain of influence (cf. [CH]).

Consider the Cauchy problem for the equation (\ref{eq:1.1}) in half-space  $x_0>t_0$  with
 initial conditions
\begin{equation}                                       \label{eq:2.1}
u(t_0,x)=\varphi_0(x),\ \ u_{x_0}(t_0,x)=\varphi_1(x),
\end{equation}
where $\mbox{supp\ }\varphi_k(x)\subset\overline{X},\ k=0,1,\ x\in\R^n$.
Denote by $D_+(\varphi_0,\varphi_1)$  the support of $u(x_0,x)$  for $x_0\geq t_0$.
Let $D_+(X\times \{x_0=t_0\})$ be the closure of the union of $D_+(\varphi_0,\varphi_1)$  over
all $\varphi_0,\varphi_1$ with supports on $\overline{X}$.  Then $D_+(X\times\{x_0=t_0\})$  is 
the forward domain of influence of $X\times\{x_0=t_0\}$.
Therefore $\Omega_{int}\times\R$  is a black hole if
$D_+(\Omega_{int}\times\{x_0=t_0\})\subset\overline{\Omega}_{int}\times \R$,  and 
$\Omega_{int}\times\R$ 
is a white hole if 
$D_+(\Omega_{ext}\times\{x_0=t_0\})\subset\overline{\Omega}_{ext}\times \R$   
  for all $t_0\in \R$.  

There is a geometric description of $D_+(X\times\{x_0=t_0\})$.

Let $[g_{jk}(x)]_{j,k=0}^n$  be the metric tensor corresponding to the operator $(1,1)$.   
Consider a curve in $\R^{n+1}$:
\begin{equation}                                              \label{eq:2.2}
x_0=x_0(s),\  x=x(s),\ s\geq 0,\ x_0(0)=y_0,\ x(0)=y.
\end{equation}
The curve (\ref{eq:2.2})  is called a time-like ray if
\begin{equation}                                                  \label{eq:2.3}
\sum_{j,k=0}^ng_{jk}(x(s))\frac{dx_j(s)}{ds}\frac{dx_k(s)}{ds}>0,
\  \ \frac{dx_0}{ds}>0,
\ \ \mbox{for\ }\ \ s\geq 0.
\end{equation}
\begin{theorem}(cf., for example, [CH]).                        \label{theo:2.1}
The forward domain of influence $D_+(X\times\{x_0=t_0\})$  is the closure of the union  
of all time-like rays 
starting at $\overline X\times\{x_0=t_0\}$.
\end{theorem}
Let $\{S_0(x)=0\}\times\R$  be a characteristic surface for (\ref{eq:1.1}),  i.e.
\begin{equation}                                                \label{eq:2.4}
\sum_{j,k=1}^n g^{jk}(x)S_{0x_j}(x) S_{0x_k}=0 \ \ \mbox{when\ } S_0(x)=0.
\end{equation}
We assume that $S_{0x}(x)$ is the outward normal to $S_0(x)=0,
\ S_{0x}(x)\neq 0$  when  $S_0(x)=0$.

\begin{theorem}(cf. [E2]).                                    \label{theo:2.2}
The domain $\Omega_{int}\times\R$  is a black hole if (\ref{eq:2.4}) holds and
\begin{equation}                                                \label{eq:2.5}
\sum_{j=1}^n g^{j0}(x)S_{0x_j}(x)<0 \ \ \mbox{when\ } S_0(x)=0,
\end{equation}
and the domain $\Omega_{int}\times\R$  is a white hole if (\ref{eq:2.4}) holds
and if
\begin{equation}                                                \label{eq:2.6}
\sum_{j=1}^n g^{j0}(x)S_{0x_j}(x)>0 \ \ \mbox{when\ } S_0(x)=0.
\end{equation}
\end{theorem}

The proof of Theorem \ref{theo:2.2} based on the study of the time-like rays starting on $S_0(x)=0$
was given in [E2].

One can also prove Theorem \ref{theo:2.2}  using the energy-type estimates of the solutions of
the equation (\ref{eq:1.1})
in the  Sobolev spaces.
\begin{theorem}(cf., for example, [E4]).                               \label{theo:2.3}
Let (\ref{eq:2.4}) and (\ref{eq:2.5}) hold.   Then for any solution  
$u(x_0,x)$  of (\ref{eq:1.1})  we have
\begin{equation}                                           \label{eq:2.7}
\|u(x_0,\cdot)\|_{1,\Omega_{ext}}^2 + 
\|u_{x_0}(x_0,\cdot)\|_{0,\Omega_{ext}}^2
\leq
C_T(\|u(t_0,\cdot)\|_{1,\Omega_{ext}}^2 + 
\|u_{x_0}(t_0,\cdot)\|_{0,\Omega_{ext}}^2),
\end{equation}
where  $T$  is arbitrary, $t_0\leq x_0\leq T,\ \|v\|_{p,\Omega_{ext}}$  is the norm
in the Sobolev space $H_p(\Omega_{ext})$.
\end{theorem}
Having the estimate (\ref{eq:2.7})  it is easy to see that  $\Omega_{int}\times\R$  is a black hole:

Consider the Cauchy problem (\ref{eq:1.1}), (\ref{eq:2.1}) in the half-space $x_0>t_0$,
where $\mbox{supp\ }\varphi_k(x)\subset\overline\Omega_{int},\ k=0,1$.  Then
$u(t_0,x)=u_{x_0}(t_0,x)=0$  when $x\in \Omega_{ext}$.  Applying the estimate (\ref{eq:2.7})
  we get that $u(x_0,x)=0$  in $\Omega_{ext}\times(t_0,+\infty)$.  
Therefore $\mbox{supp\ }u(x_0,x)\subset \overline{\Omega}_{int}\times[t_0,+\infty)$,  i.e. 
$D_+(\Omega_{int}\times\{x_0=t_0\})\subset\overline{\Omega}_{int}\times\R$  for any $t_0$,
i.e.  $\Omega_{int}\times\R$  is a black hole.

Analogous result holds for the white hole.
\begin{theorem}(cf., for example, [E4]).                               \label{theo:2.4}
Let (\ref{eq:2.4}) and (\ref{eq:2.6}) hold.   Then for any solution  
$u(x_0,x)$  of (\ref{eq:1.1})  we have
\begin{equation}                                           \label{eq:2.8}
\|u(x_0,\cdot)\|_{1,\Omega_{int}}^2 + 
\|u_{x_0}(x_0,\cdot)\|_{0,\Omega_{int}}^2
\leq
C_T(\|u(t_0,\cdot)\|_{1,\Omega_{int}}^2 + 
\|u_{x_0}(t_0,\cdot)\|_{0,\Omega_{int}}^2),
\end{equation}
where  $T$  is arbitrary, $t_0\leq x_0\leq T$.
\end{theorem}
As in the case of Theorem \ref{theo:2.3}  the estimate (\ref{eq:2.8})  implies that 
$D_+(\Omega_{ext}\times\{x_0=t_0\})\subset\overline{\Omega}_{ext}\times\R$,
i.e. $\Omega_{int}\times\R$  is a white hole.

\section{The ergosphere.}
\label{section 3}
\init

The ergosphere is the surface $S(x)=0$  where
\begin{equation}                                   \label{eq:3.1}
g_{00}(x)=0.
\end{equation}
We assume that $S(x)=0$  is a closed smooth surface,  
$g_{00}(x)>0$  in the exterior of $S(x)=0$  and $g_{00}(x)<0$  in the interior 
of $S(x)=0$  near $S(x)=0$.  
We say  that $S(x)=0$  is a smooth  surface  if  $S_x(x)\neq 0$  when  $S(x)=0$.

Let $\Delta(x)=\det[g^{jk}(x)]_{j,k=1}^n$.  
It is easy to show  (cf. [E1])
that $g_{00}(x)=0$  if and only if $\Delta(x)=0$.  In the case of the Gordon equation  
(cf. (\ref{eq:1.4}))
the equation of the ergosphere 
is 
\begin{equation}                                       \label{eq:3.2}
|w(x)|^2=\frac{c^2}{n^2(x)}.
\end{equation}
Note that $(0,\xi)$  is not a characteristic direction for (\ref{eq:1.1}) 
for any $\xi\neq 0$  when $g_{00}(x)>0$.
It may happen that the ergosphere $S(x)=0$  is also a characteristic surface,  i.e. 
$\{S(x)=0\}\times\R$  is an event horizon.  The celebrated example of such situation is 
the Schwarzschield black hole.  The Schwarzchield  metric has the following form in 
the Cartesian coordinates (cf. [V2]):
\begin{equation}                                        \label{eq:3.3}
ds^2=(1-\frac{2m}{R})dt^2-dx^2-dy^2-dz^2-\frac{4m}{R}dtdR-\frac{2m}{R}(dR)^2,
\end{equation}
where $R=\sqrt{x^2+y^2+z^2}$.  Therefore
\begin{equation}                                           \label{eq:3.4}
g_{00}=1-\frac{2m}{R}=0
\end{equation}
is the ergosphere.  Note that $R=2m$  is also  a characteristic surface and $\{R=2m\}\times\R$
is a black hole.  We shall call the black hole such that the ergosphere is also an event horizon 
the Schwarzschield type black holes.

Note that Schwarzschield  type black holes is unstable with respect to the perturbations of
metrics (see \S 5 below).  If we perturb the metric the ergosphere persists since we assume
that the equation (\ref{eq:3.1}) is smooth,  i.e.  it has a non-vanishing normal at any point.
However the perturbed surface may cease to be a characteristic  
surface and there is no characteristic surface near by (cf. \S 5).

It is easier to study the behavior of the solutions of (\ref{eq:1.1}) in the exterior of
the Schwarzschield type black hole 
than in the case when the black hole is inside the ergosphere.  For example,
an important problem in the general relativity studied in [DR]  is the uniform boundedness 
of solutions of (\ref{eq:1.1}) in the exterior of the black hole.  Using the ideas from
[DR]  one can prove the following theorem:
\begin{theorem}(cf. [E3])                                         \label{theo:3.1}
Let $\{S_0(x)=0\}\times\R$  be the ergosphere and the boundary of a black hole,
$S_{0x}(x)\neq 0$  when  $S_0(x)=0$.
Consider the Cauchy problem for (\ref{eq:1.1})  in the exterior $\Omega_{ext}\times(0,+\infty)$  of
the black hole with the initial conditions
$$
u(0,x)=\varphi_0,\ \ u_{x_0}(0,x)=\varphi_1(x),\ \ x\in\Omega_{ext},
$$
where $\varphi_0(x),\varphi_1(x)$  are smooth and rapidly 
decaying when $|x|\rightarrow\infty$.
Then $u(x_0,x)$  is uniformly bounded in $\Omega_{ext}\times(0,+\infty)$:
\begin{equation}                                             \label{eq:3.5}
|u(x_0,x)|\leq C.
\end{equation}
\end{theorem}

 \section{The case of two space dimensions.}
\label{section 4}
\init         

Let $S$  be an ergosphere,  i.e. $\Delta(x)=g^{11}(x)g^{22}(x)-(g^{12}(x))^2=0, \ x=(x_1,x_2)\in \R^2$.  
We assume that $S$  is a smooth Jordan curve.  Let $S_1$  be another closed Jordan curve 
inside $S$.  Denote by $\Omega$  the region between $S$  and $S_1$.  We assume that 
$\Delta(x)<0$  in $\overline{\Omega}\setminus S$.  Let $K_+(y)$  be the half-cone of all forward 
time-like 
directions at $y\in S_1$,  i.e. $K_+(y)=\{(\dot{x}_0,\dot{x}_1,\dot{x}_2)\in \R^3:
\sum_{j,k=0}^2 g_{jk}(y)\dot{x}_j\dot{x}_k>0,\dot{x}_0>0\}$.  Let
$N(y)$  be the outward normal to $S_1$.  We assume that either
\begin{eqnarray}                               \label{eq:4.1}
(\dot{x}_0,\dot{x}_1,\dot{x}_2)\cdot (0,N(y))>0
\\
\mbox{or}
\nonumber
\\
\nonumber
(\dot{x}_0,\dot{x}_1,\dot{x}_2)\cdot (0,N(y))<0
\end{eqnarray}for all $y\in S_1$  and all $(\dot{x}_0,\dot{x}_1,\dot{x}_2)\in \overline K_+(y)$.

{\bf Remark 4.1}
The interior of $S_1$  is called a trapped region if \\
$(\dot{x}_0,\dot{x}_1,\dot{x}_2)\cdot(0,N(y))<0$
for all $y\in S_1$  and $(\dot{x}_0,\dot{x}_1,\dot{x}_2)\in \overline{K}_+(y)$.

The main result of [E2]  is the following theorem:
\begin{theorem}                                                \label{theo:4.1}
Suppose the ergosphere $S$  is not characteristic for any $x\in S$,  i.e.
\begin{equation}                                             \label{eq:4.2}
\sum_{j,k=1}^2g^{jk}(x)\nu_j(x)\nu_k(x)\neq 0,\ \ \forall x\in S,
\end{equation} 
where $\nu(x)=(\nu_1(x),\nu_2(x))$  is the normal to $S$ at $x\in S$.
Suppose 
the condition (\ref{eq:4.1})  is satisfied on $S_1$.  Then there exists a smooth Jordan curve 
$S_0(x)=0$  between $S$  and $S_1$  such that $\{S_0(x)=0\}\times \R$  is a characteristic surface,
i.e. $\{S_0(x)=0\}\times R$  is an event horizon.
\end{theorem}

\underline{Sketch of the proof of Theorem \ref{theo:4.1}}.
\\
\\
Since $\Delta(x)<0$  in $\Omega$  there exist two families 
$S^\pm(x)=\mbox{const} $ 
of characteristic curves 
$$
\sum_{j,k=1}^2g^{jk}(x)S_{x_j}^\pm(x)S_{x_k}^\pm(x)=0.
$$
One can construct two vector fields $f^\pm(x)=(f_1^\pm(x),f_2^\pm(x))$  on $\overline{\Omega}$
such that $f^\pm(x)\neq (0,0)$  for all $x\in \overline{\Omega},\ f^+(x)\neq f^-(x)$
on $\overline{\Omega}\setminus S,\ f^+(x)= f^-(x)$  on $S$ and
\begin{equation}                               \label{eq:4.3}
f_1^\pm(x)S_{x_1}^\pm (x)+f_2^\pm(x)S_{x_2}^\pm(x)=0 \ \ \mbox{in}\ \ \overline{\Omega}.
\end{equation}
Consider two systems of differential equations:
\begin{equation}                                  \label{eq:4.4}
\frac{dx^\pm(\sigma)}{d\sigma}=f^\pm (x(\sigma)),\ \ \sigma\geq 0,\ \ x^\pm(0)=y\in S.
\end{equation}
Let
\begin{eqnarray}                                     \label{eq:4.5}
\frac{dx_j(s)}{ds}=2\sum_{k=0}^2 g^{jk}(x(s))\xi_k(s),\ \ x_j(0)=y_j,\ \ 0\leq j\leq 2,
\\
\nonumber
 \frac{d\xi_p(s)}{ds}=-\sum_{j,k=0}^2 g_{x_p}^{jk}(x(s))\xi_j(s)\xi_k(s),\ \ \xi_p(0)=\eta_p,
\ \ 0\leq p\leq 2,
\end{eqnarray}
be the equations of null-bicharacteristics for (\ref{eq:1.1}).                              
Note that $\frac{d\xi_0(s)}{ds}=0$,  i.e.  $\xi_0(s)=\eta_0$  for all $s\geq 0$  
and we take $\eta_0=0$.
Therefore (\ref{eq:4.5})  is a null-bicharacteristic  if $\sum_{j,k=1}^2 g^{jk}(y)\eta_j\eta_k=0$.

It can be shown that the curves $x=x^\pm(\sigma)$  of (\ref{eq:4.4}) are the projections on 
$(x_1,x_2)$-plane of some forward null-bicharacteristics.  Since $\frac{dx_0(s)}{ds}\neq 0$  on 
these bicharacteristics we can use $x_0$  as a parameter instead of $\sigma$.  It appears that 
for one family (say $x=x^+(\sigma)$)  $\sigma$ is decreasing when $x_0$  is increasing and
for another family ($x=x^-(\sigma)$)  $\sigma$  is increasing when $x_0$  is increasing.

The condition (\ref{eq:4.1})  
is equivalent to the condition that the projections on
$(x_1,x_2)$-plane of all forward null-bicharacteristics are either leaving $\Omega$  when $x_0$ is 
increasing or are entering $\Omega$ when $x_0$  is increasing.
Suppose for the definiteness that the projections of all null-bicharacteristics are leaving $\Omega$  
when $x_0$  is increasing.  Then the trajectory $x=x^-(\sigma)$  starting on $S$ can not
 reach $S_1$ (cf. [E2]).
Therefore by the Poincare-Bendixson theorem there exists a limit cycle $S_0(x)=0$,  i.e.
a closed Jordan curve in $\Omega$  that is a characteristic curve.  Therefore 
$\{S_0(x)=0\}\times\R$  is an event horizon,  i.e. a boundary of either black or white hole.

{\bf Remark 4.2}
Since conditions (\ref{eq:4.1}),  (\ref{eq:4.2})  hold  when we slightly perturb the metric
the black and white holes obtained by Theorem \ref{theo:4.1} are stable.

{\bf Example 4.1}  (cf.  [V1])  Consider the acoustic equation with 
the metric (\ref{eq:1.5})
when $n=2,\ \rho=c=1$,
\begin{equation}                                   \label{eq:4.6}
v(x)=(v^1(x),v^2(x))=\frac{A}{r}\hat r+\frac{B}{r}\hat\theta,
\end{equation}
where $r=|x|=\sqrt{x_1^2+x_2^2},\ \hat r=\frac{x}{|x|},\ \hat\theta=(-\frac{x_2}{|x|},\ \frac{x_1}{|x|}),\
A>0,\ B>0$  are constants.  We assume that $\sqrt{A^2+B^2}\geq |x|>r_1$,  where $r_1<A$.  
The ergosphere in this case is $r=\sqrt{A^2+B^2}$.  The differential equations (\ref{eq:4.4}) have
the following form in the polar coordinates $(r,\theta)$:  
\begin{equation}                                          \label{eq:4.7}
\frac{dr}{ds}=A^2-r^2,\ \ \frac{d\theta}{ds}=\frac{AB}{r}+\sqrt{A^2+B^2-r^2},
\end{equation}
and
\begin{equation}                                          \label{eq:4.8}
\frac{dr}{ds}=-1,\ \ \frac{d\theta}{ds}=\frac{1-\frac{B^2}{r^2}}{\frac{AB}{r}+\sqrt{A^2+B^2-r^2}}.
\end{equation}
It follows from (\ref{eq:4.7}) that $r=A$  is a limit cycle and $\{r<A\}\times\R$  is a white hole.

{\bf Example 4.2}  Consider the same situation as in Example 4.1 
with $B=0$  and the domain $r<A$.  Then $r=A$  is the ergosphere and $\{r=A\}\times\R$  is the event horizon.
Since $A>0$  we have that $\{r<A\}\times \R$ is a white hole.  Note that the equations (\ref{eq:4.4}) have
the following form in the polar coordinates:
\begin{equation}                                        \label{eq:4.9}
\frac{dr}{d\theta}=\pm\sqrt{A^2-r^2},\ \ r(\theta_0)=A.
\end{equation}
It has a solution $r=A$,  which is   the event horizon,  and it  also has other solutions 
$r=A\cos(\theta-\theta_0)$  that touch the event horizon at $\theta=\theta_0$.  In 
general situation when the  
ergosphere coincides with the event  
horizon  the solution of (\ref{eq:4.4}) are also tangent to the event horizon.

{\bf Remark 4.3}
When $B\neq 0$  is small 
Example 4.1 can be viewed as a perturbation of Example 4.2.  Note that when $B$  is small 
the stable event horizon $r=A$  will be close to the ergosphere $r=\sqrt{A^2+B^2}$.
 A similar situation will happen in the general case: 

If the left hand side of the 
condition (\ref{eq:4.2}) is small then the stable event horizon obtained in Theorem \ref{theo:4.1}
will be close to the ergosphere $S$.  If the left hand side of (\ref{eq:4.2}) changes sign on $S$ 
(for example, if $B$ in (\ref{eq:4.6}) depends on $\theta$  and changes sign when
$0\leq \theta\leq 2\pi$)  there will be no event horizon near $S$.

\section{Axisymmetric metrics and rotating black holes.}
\label{section 5}
\init

Let $(\rho,\varphi,z)$  be the cylindrical coordinates in $\R^3$:
\begin{equation}                                \label{eq:5.1}
x=\rho\cos\varphi,\ \ y=\rho\sin\varphi,\ \ z=z.
\end{equation}
A stationary axisymmetric metric in $\R^3\times\R$  is the metric that does not depend on 
$t$  and $\varphi$.
For the convenience,  we shall use the following notations:
\begin{equation}                                  \label{eq:5.2}
y_0=t,\ \ y_1=\rho,\ \ y_2=z,\ \ y_3=\varphi.
\end{equation}
Then the stationary axisymmetric metric has the form:
\begin{equation}                                  \label{eq:5.3}
ds^2=\sum_{j,k=0}^3g_{jk}(\rho,z)dy_jdy_k,
\end{equation} 
where $g_{jk}(\rho,z)$  are smooth and even in $\rho$.  Denote
$$
[g^{jk}(\rho,z)]_{j,k=0}^3=\left([g_{jk}(\rho,z)]_{j,k=0}^3\right)^{-1}.
$$
The ergosphere is given by the equation
\begin{equation}                                 \label{eq:5.4}
g_{00}(\rho,z)=0,
\end{equation}
or,  equivalently:
\begin{equation}                                     \label{eq:5.5}
\Delta(\rho,z)=\det [g^{jk}(\rho,z)]_{j,k=1}^3=0.
\end{equation}

We will be looking for the rotating black and white holes,  i.e.  when the event horizon 
has the form:
\begin{equation}                                       \label{eq:5.6}
\{S(\rho,z)=0\}\times S^1\times\R,
\end{equation}
where  $S(\rho,z)=0$  is a closed smooth curve in the $(\rho,z)$-plane,  even in $\rho,
\ \varphi\in S^1,\ t\in \R,\ S^1$  is the unit circle.  More precisely,  we have to take
in (\ref{eq:5.5})  the restriction of the curve $S(\rho,z)=0$  to the half-plane $\rho\geq 0$
but we did not indicate this in (\ref{eq:5.6})  for the simplicity of notation.  Since
(\ref{eq:5.6})  is the event horizon it must be a characteristic surface,  i.e.
\begin{equation}                                     \label{eq:5.7}
\sum_{j,k=1}^2g^{jk}(\rho,z)S_{y_j}(\rho,z)S_{y_k}(\rho,z)=0\ \ \ \mbox{on}\ \ S(\rho,z)=0.
\end{equation}
Here $y_1=\rho,\ y_2=z$.  Therefore 
$\{S(\rho,z)=0\}\times\R$  is  the event horizon for the tensor 
$[g^{jk}(\rho,z)]_{j,k=0}^2$,  i.e.  for the case of two dimensions considered in the previous section.

Define 
\begin{equation}                                         \label{eq:5.8}
\Delta_1(\rho,z)=\det [g^{jk}(\rho,z)]_{j,k=1}^2
=g^{11}(\rho,z)g^{22}(\rho,z)-(g^{12}(\rho,z))^2.
\end{equation}
We shall call the curve $\Delta_1(\rho,z)=0$  the restricted ergosphere since it is the ergosphere 
of the two-dimensional 
problem for $[g^{jk}(\rho,z)]_{j,k=0}^2$.  We can extend all results of \S 4 to the case of rotating 
black and white holes.  For example,
let the curve $\Delta_1(\rho,z)=0$  be a Jordan curve  such that 
$$
\sum_{j,k=1}^2 g^{jk}(\rho,z)\nu_j(\rho,z)\nu_k(\rho,z)\neq 0 \ \ \mbox{on}\ \ \Delta_1(\rho,z)=0,
$$
and let $S_1$  be a Jordan curve inside $\Delta_1(\rho,z)=0$.  

Suppose  conditions (\ref{eq:4.1})  are satisfied where the matrix $[\tilde g_{jk}(\rho,z)]_{j,k=0}^2$ 
is the inverse to $[g^{jk}(\rho,z)]_{j,k=0}^2$.  Then there exists a Jordan curve $S_0(\rho,z)=0$
between $\Delta_1=0$  and $S_1$  such that $\{S_0(\rho,z)=0\}\times S^1\times \R$ is the event 
horizon in $\R^3\times \R$.

Consider now the problem of the stability of the black and white holes with respect to 
the perturbations of metrics.

The famous example of an axisymmetric metric is the Kerr metric. 
The Kerr metric  in the Kerr-Schild coordinates has the  form (see [V2]):
\begin{eqnarray}                                          \label{eq:5.9}
ds^2=dt^2-dx^2-dy^2-dz^2 \ \ \ \ \ \ \ \ \ \ \ \ \ \ \ \ \
\\
\nonumber
-\ \ \frac{2mr^3}{r^4+a^2z^2}
\left[dt+\frac{r(xdx+ydy)}{r^2+a^2}+
\frac{a(ydx-xdy)}{r^2+a^2}
+\frac{z}{r}dz\right]^2,
\end{eqnarray}
where
\begin{equation}                               \label{eq:5.10}
r(x,y,z)=\sqrt{\frac{(R^2-a^2)+\sqrt{(R^2-a^2)^2+4a^2z^2}}{2}},\ \ R^2=x^2+y^2+z^2.
\end{equation}
It follows from (\ref{eq:5.4})  and (\ref{eq:5.9}) that the ergosphere is
\begin{equation}                               \label{eq:5.11}
r^4+a^2z^2-2mr^3=0.
\end{equation}
One  can  show   that (\ref{eq:5.11}) consists of two curves  in $(\rho,z)$-plane
\begin{equation}                                \label{eq:5.12}
r-\left(m+\sqrt{m^2-\frac{a^2z^2}{r^2}}\right)=0,
\end{equation}
\begin{equation}                                \label{eq:5.13}
r-\left(m-\sqrt{m^2-\frac{a^2z^2}{r^2}}\right)=0,
\end{equation}
Equation (\ref{eq:5.12})  defines the outer ergosphere and (\ref{eq:5.13})  defines the inner 
ergosphere  for the Kerr metric.

Compute the restricted ergosphere $\Delta_1(\rho,z)$ (cf. (\ref{eq:5.8})) for the Kerr metric.
The inverse to the Kerr metric tensor has the form: 
\begin{equation}                              \label{eq:5.14}
\eta^{jk}+\frac{2mr^3}{r^4+a^2z^2}l^jl^k,
\end{equation}
where
$$
(l^0,l^1,l^2,l^3)=\left(-1,\frac{rx+ay}{r^2+a^2},\frac{ry-ax}{r^2+a^2},\frac{z}{r}\right).
$$
 In  the $(\rho,z,\varphi)$  coordinates we have
\begin{equation}                              \label{eq:5.15}
g^{jk}(\rho,z)=\xi^{jk}+\frac{2mr^3}{r^4+a^2z^2}m^jm^k,
\end{equation}
where  $(m^0,m^1,m^2,m^3)=\left(-1,\frac{r\rho}{r^2+a^2},\frac{z}{r},\frac{-a}{r^2+a^2}\right),\
\xi^{jk}$  is $\eta^{jk}$  in the cylindrical coordinates, 
$\xi^{00}=1,\ \xi^{11}=\xi^{22}=-1,\ \xi^{33}=-\frac{1}{\rho^2},\ 
\xi^{jk}=0$  for $j\neq k$.  
Therefore
\begin{equation}                                  \label{eq:5.16}
\Delta_1(\rho,z)=1-\frac{2mr^5\rho^2}{(r^4+a^2z^2)(r^2+a^2)^2}
-\frac{2mrz^2}{r^4+a^2z^2}.
\end{equation}
The equation $\Delta_1=0$  for the Kerr metric can be substantially simplified.

\begin{pro}(cf. [E3])                                       \label{prop:5.1}
The equation $\Delta_1(\rho,z)=0$  is equivalent 
to two equations $r-r_+=0$ and $r-r_-=0$  where 
$r_\pm=m\pm\sqrt{m^2-a^2}$.
\end{pro}
It happens that $r=r_\pm$  are two event horizons,  $r-r_+=0$  is called the outer event horizon and
$r-r_-=0$  is  the inner event horizon.  More exactly, $\{r=r_\pm\}\times S^1\times \R$ 
 are the event horizons.  
Therefore $\Delta_1(\rho,z)=0$  relates explicitly the event horizons of the Kerr metric to 
the metric tensor.

{\bf Definition 5.1}
Let $\{\psi=0\}\times S^1\times\R$  be the event horizon for the metric $[g_{jk}]_{j,k=0}^3$.  
We say that this event horizon is stable in the class of axisymmetric metrics if any smooth
family $[g_{\e jk}]$  of axisymmetric metrics, $0\leq \e\leq\e_0,\ g_{0jk}=g_{jk}$,  has a smooth 
family of event horizons $\{\psi_\e=0\}\times S^1\times\R$  such that
$\psi_0=\psi$.  Otherwise we say that $\{\psi=0\}\times S^1\times\R$  is 
an unstable event horizon.

We restrict perturbations to a more narrow class of  axisymmetric metrics of the form:
\begin{equation}                                   \label{eq:5.17}
g_\e^{jk}=\xi^{jk}+v_\e^j(\rho,z)v_\e^k(\rho,z).
\end{equation}
Note that the Kerr metric and the metrics (\ref{eq:1.4}), (\ref{eq:1.5})  
have the form (\ref{eq:5.17}).

\begin{pro}(cf. [E3])                                          \label{prop:5.2}
Let $\Delta_1=0$  be the restricted ergosphere,  and let $\Delta_1=0$  
be a characteristic  curve, i.e.  
$\{\Delta_1=0\}\times S^1\times \R$  is an event horizon.  
Then this event horizon is unstable when we consider perturbations 
in the class of the metrics of the form (\ref{eq:5.17}).  
In particular,  the outer and the inner event horizons for the Kerr metric are unstable.                                    
\end{pro}

\underline{Sketch of the proof of Proposition \ref{prop:5.2}.}
\\
\\
Denote $\Delta_1^\e=g_\e^{11}g_\e^{22}-(g_\e^{12})^2$.  Then $\Delta_1^\e =0$  is a smooth perturbation 
of restricted ergosphere,  $\Delta_1^0=\Delta_1$.   We can choose perturbations of the form  
 (\ref{eq:5.17})  such that $\Delta_1^\e=0$  will not be a characteristic curve  for $0< \e\leq \e_0$.
Moreover,  one
can choose $[g_{\e jk}]$  such that there is no characteristic curve near $\Delta_1^\e=0$ (cf. [E3]).

In the next proposition we shall prove that there  is a rich class of perturbations 
$[g_{\e jk}],\ 0\leq\e\leq \e_0,$  of the Kerr metric that have a smooth family  of event horizons
$\{\Delta_1^\e=0\}\times S^1\times\R$  such that $\{\Delta_1^0=0\}\times S^1\times \R$  is the Kerr 
event horizon.  

\begin{pro}(cf. [E3])                                    \label{prop:5.3}
Let $\Delta_1=0$  be a restricted ergosphere and $\{\Delta_1=0\}\times S^1\times \R$  is an event horizon.  
Let $\Delta_1^\e=0$  be arbitrary family of closed, even in $\rho$,  smooth curves,  $0\leq \e\leq \e_0$,  
such  that $\Delta_1^0=\Delta_1$.
Then there exists a family of metrics $[g_{\e jk}],\ 0\leq \e\leq \e_0,\ [g_{0jk}]=[g_{jk}]$,  
of the form  (\ref{eq:5.17})  such that $\Delta_1^\e=0$  are restricted ergospheres for $[g_{\e jk}]$  and
$\{\Delta_1^\e=0\}\times S^1\times\R$ are event horizons for $[g_{\e jk}]$.
\end{pro}

{\bf Remark 5.1}
We shall call the event horizons obtained in Theorem \ref{theo:4.1}  the stable event horizons and 
the event horizons 
that coincide with $\{\Delta_1=0\}\times S^1\times\R$  the Schwarzschield type event horizons.
According to this definition the Kerr event horizon is a Schwarzschield type event horizon.

Fix some axisymmetric 
metric $[g_{jk}]$ that has a Schwarzschield type event horizon $\{\Delta_1= 0\}\times S^1\times\R$.
The proposition \ref{prop:5.3}  shows that in any neighborhood of $[g_{jk}]$  there are metrics having 
Schwarzschield type event horizons close  to $\{\Delta_1=0\}\times S^1\times \R$.  Also it follows
from the Remark 4.3  that in any neighborhood of $[g_{jk}]$  (in particular,  in any neighborhood of
the Kerr metric)  there are stable event horizons close to $\{\Delta_1=0\}\times S^1\times \R$.

\section{Determination of the ergosphere by the boundary measurements.}
\label{section 6}
\init

Let $u(x_0,x)$  be the solution of (\ref{eq:1.1})  in a cylinder $\Omega\times \R$  satisfying 
the zero initial 
conditions
\begin{equation}                                              \label{eq:6.1} 
u=0\ \ \ \mbox{for}\ \ x_0\ll 0,\ \ x\in \Omega,
\end{equation}
and the boundary condition
\begin{equation}                                            \label{eq:6.2}
u|_{\partial\Omega\times\R}=f.
\end{equation}
Here $\Omega$  is a smooth bounded domain in $\R^n, \ f$  is a smooth function with a compact 
support in $\partial\Omega\times\R$.  The solution of the initial-boundary problem (\ref{eq:1.1}),
(\ref{eq:6.1}),  (\ref{eq:6.2})  exists and is unique assuming that $\partial\Omega\times\R$  is not 
characteristic at any point  and $g^{00}(x)>0$  on $\overline{\Omega}$.

Denote by $\Lambda f$  the following operator (the DN operator):
\begin{equation}                                    \label{eq:6.3}
\Lambda f=\sum_{j,k=1}^ng^{jk}(x)\frac{\partial u}{\partial x_j}\nu_k(x)\left.\left(\sum_{p,r=1}^n
g^{pr}(x)\nu_p\nu_r\right)^{-\frac{1}{2}}\right|_{\partial\Omega\times\R},
\end{equation}
where
$(\nu_1(x),...,\nu_n(x))$
is the  outward unit normal to $\partial\Omega$.  
Let $\Gamma$  be an open subset of $\partial\Omega$.  We say that boundary measurements on
$\Gamma\times(0,T)$  are taken if we know $\Lambda f$  on $\Gamma\times(0,T)$  for all
$f$ with support in  $\overline{\Gamma}\times[0,T]$.

The inverse boundary value problem is the determination of $[g^{jk}(x)]_{j,k=0}^n$  knowing
the boundary measurements on $\Gamma\times(0,T).$
Let 
\begin{equation}                                          \label{eq:6.4}
\hat x=\varphi(x),\ \ \hat x_0=x_0+a(x),
\end{equation}
where  $\hat x=\varphi(x)$  is a diffeomorphism of $\overline\Omega$  onto a new domain $\overline\Omega,
\ a(x)\in C^\infty(\overline\Omega)$.  We assume that
\begin{equation}                                          \label{eq:6.5}
\varphi(x)=x \ \ \mbox{on}\ \ \overline\Gamma,\ \ a(x)=0 \ \ \mbox{on} \ \ \overline\Gamma.
\end{equation}
Note that (\ref{eq:6.4})  transform (\ref{eq:1.1})  to an equation of the same form with
a new tensor $[\hat g^{jk}(\hat x)]_{j,k=0}^n$  isometric to the old one.   It follows from
(\ref{eq:6.5})  that DN operator $\Lambda$  does not change under the change of variables
(\ref{eq:6.4}),  (\ref{eq:6.5}).  If there exists an event horizon inside $\Omega\times\R$  then
we can not determine the metric inside the event horizon  since any change of metric inside 
the event horizon will not change boundary measurements. But we can try to recover the event
horizon itself (up to diffeomorphism (\ref{eq:6.4}),  (\ref{eq:6.5})).

This is an open problem.  We can prove only  that the boundary measurements allow to determine
the ergosphere.

\begin{theorem}(cf. [E3])                                               \label{theo:6.1}
Consider the wave equation (\ref{eq:1.1}).  Assume that $g^{00}(x)>0$  on 
$\overline\Omega$  and the normal to $\partial\Omega$  is not characteristic 
at any point of $\partial\Omega$.  
Let $\Delta(x)=0$  be the ergosphere,  $\Delta(x)=0$  is a smooth closed surface,
$\Delta(x)>0$  in $\overline{\Omega}$  outside of $\Delta(x)=0$.  Let $\Gamma$  be
an open subset of $\partial\Omega$.  Then the boundary measurements on $\Gamma\times(0,+\infty)$
determine $\Delta(x)=0$ up to the change of variables (\ref{eq:6.4}),  (\ref{eq:6.5}).
\end{theorem} 
Note that for the proof of Theorem \ref{theo:6.1} it does not matter 
whether the ergosphere is an event horizon or not.  The proof 
is an extension of the proof of Theorem 2.3 in [E1].  

 We will determine the ergosphere  by determining
 the metric in 
$\overline\Omega\cap\overline\Omega_{ext}$,   where $\overline\Omega_{ext}$ is
the exterior  of $\Delta(x)=0$.

We start with the determination of the metric in a small neighborhood of $\Gamma$ and
gradually continue to recover the metric deeper in $\Omega$.  As we progress the time interval
$(0,T)$  needed to reach the point $x\in\Omega$  increases when the point approaches the ergosphere.
One can show that $T\rightarrow +\infty$  when $x\rightarrow\{\Delta(x)=0\}$.  This is the reason why
one needs the unlimited time interval $(0,+\infty)$  to recover the ergosphere. 

{\bf Remark 6.1}  Let $L(x,\frac{\partial}{\partial x},\frac{\partial}{\partial x_0})u(x_0,x)=0$
be the equation (\ref{eq:1.1})  in $\R^{n+1}$.  Making the Fourier transform in $x_0$  we get
\begin{equation}                                       \label{eq:6.6}
L(x,\frac{\partial}{\partial x},ik)\tilde u(k,x_0)=0,\ \ x\in \R^n.
\end{equation}
Suppose that 
\begin{equation}                                       \label{eq:6.7}
L(x,\frac{\partial}{\partial x},ik)=-\Delta-k^2\ \ \ \mbox{for}\ \ |x|>R.
\end{equation}
Let $a(\theta,\omega,k)$  be the scattering
amplitude for the operator (\ref{eq:6.6}).  It is well known that the scattering amplitude given for all 
$k>0,\ \theta\in S^{n-1},\ \omega\in S^{n-1},$ 
determines the DN operator on the $\{|x|=R\}\times[0,+\infty)$.  Therefore by the Theorem \ref{theo:6.1}
the scattering amplitude determines the ergosphere. Note that when (\ref{eq:6.7}) holds $a(\theta,\omega,k)$
is real analytic in $(\theta,\omega,k)$.  Therefore it is enough to know $a(\theta,\omega,k)$ in
an arbitrary neighborhood of some point $(\theta_0,\omega_0,k_0)$  to determine the ergosphere.


\begin{thebibliography}{9999}
\bibitem[C]{} Chandrasekhar, S., The Mathematical Theory of Black Holes,
Oxford University Press, 1998
\bibitem[CH]{} Courant, R., Hilbert, D., Methods of Mathematical Physics,
vol. II (1962),  New York, London 
\bibitem[DR]{} Dafermos, M., Rodnianski, I.,
Lectures on black holes and linear waves, arXiv:0811.0354 
\bibitem[E1]{} Eskin, G., 
Optical Aharonov-Bohm effect:  inverse hyperbolic problem approach,
Comm. Math. Phys. 284 (2008),  no 2, 317-343
\bibitem[E2]{} Eskin, G., Inverse hyperbolic problems  and optical black holes, 
Commun. Math. Phys. 297,  817-839 (2010)
\bibitem[E3]{} Eskin, G.,  Perturbations of the Kerr black hole 
and the boundness of linear waves, Journ. of Math. Phys. 51, 112501 (2010)
 \bibitem[E4]{} Eskin, G.,  Uniqueness and non-uniqueness in Inverse Hyperbolic Problems
and the Black Hole Phenomenon, Around the research of Vladimir Maz'ya  III,
Analysis and Applications - Ari Laptev (Ed.),  International Mathematical Series, 
Vol. 13, 2010, 77-104.  
\bibitem[G]{} Gordon, W., 1923, Ann. Phys.  (Leipzig) 72, 421
\bibitem[K]{} Kerr, R., Gravitational field of a spinning mass as an example of
algebraically special metrics,  Phys. Rev. Lett. 11, 237-238 (1963)
 \bibitem[LP]{} Leonhardt, U., Piwnicki, P., 1999,  Phys.  Rev. A60,  4301
\bibitem[NVV]{} Novello, M., Visser, M.,  Volovik, G. (editors),
Artificial black holes, 2002,  World Scientific,  Singapore.
\bibitem[U]{} Unruh, W., 1981, Phys. Rv. Lett., 46, 1351
\bibitem[V1]{} Visser, M., 1998,  Acoustic black holes, horizons,  ergospheres and Hawking
radiation,  Classical quantum gravity 15, No. 6, 1767-1791. 
\bibitem[V2]{} Visser, M. The Kerr space-time:  A brief introduction, arXiv:0706.0622






\end{thebibliography}
\end{document}